\documentclass[12pt]{article}

\usepackage{amsmath}
\usepackage{amssymb}
\usepackage{amsthm}
\usepackage{psfig}
\newtheorem{thm}{Theorem}[section]
\newtheorem{co}[thm]{Corollary}
\newtheorem{lem}[thm]{Lemma}

\newtheorem{assumption}[thm]{Assumption}

\newtheorem{definition}[thm]{Definition}
\newenvironment{de}{\begin{definition}\rm}{\end{definition}}
\newtheorem{example}[thm]{Example}
\newenvironment{exmp}{\begin{example}\rm}{\end{example}}
\newtheorem{remark}[thm]{Remark}
\newenvironment{rem}{\begin{remark}\rm}{\end{remark}}
\newtheorem{tab}{Table}
\newenvironment{ta}{\begin{tab}\rm}{\end{tab}}

\newcommand{\vier}[4]{\left( \begin{array}{ccc}
                   #1 &\;& #2 \\ #3 &\;& #4 \end{array} \right)}

\newcommand{\Section}[1]{\section{#1}\setcounter{equation}{0}}

\newcommand{\eqr}[1]{~\mbox{$(${\rm \ref{#1}}$)$}}
\newcommand{\tr}{{\rm tr}\,}
\newcommand{\diag}{{\rm diag}\,}

\newcommand{\V}{\mathcal{V}}
\renewcommand{\S}{\mathcal{S}}
\newcommand{\R}{\mathbb{R}}
\newcommand{\C}{\mathbb{C}}

\title{A Numerical Approach for Designing Unitary Space Time
  Codes with Large Diversity\footnote{Both authors were supported
    in part by NSF grants DMS-00-72383 and CCR-02-05310. The
    first author was also supported by a fellowship from the
    Center of Applied Mathematics at the University of Notre
    Dame. A preliminary version of this paper was presented at
    40-th Allerton Conference on Communication, Control, and
    Computing, Monticello, Illinois, October 2002.}}

\date{}

\author{Guangyue Han,\ \ Joachim Rosenthal\\
  {\normalsize Department of Mathematics}\vspace{-1mm} \\
  {\normalsize University of Notre Dame}\vspace{-1mm} \\
  {\normalsize Notre Dame, IN 46556.}\\
  {\normalsize {\em e-mail:\/} Han.13@nd.edu,\  Rosenthal.1@nd.edu}\vspace{-1mm}\\
  {\normalsize {\em URL:} http://www.nd.edu/\~{}eecoding/} }

\begin{document}\maketitle\thispagestyle{empty}

\begin{abstract}
  Unitary space-time modulation using multiple antennas promises
  reliable communication at high transmission rates. The basic
  principles are well understood and certain criteria for
  designing good unitary constellations have been presented.

  There exist two important design criteria for unitary space
  time codes. In the situation where the signal to noise ratio is
  large it is well known that the {\em diversity product} of a
  constellation should be as large as possible. It is less known
  that the {\em diversity sum} is a very important design
  criterion for codes working in a low SNR environment. For some
  special situations, it will be more practical and reasonable to
  consider a constellation optimized at certain SNR interval. For
  this reason we introduce the {\em diversity function} as a
  general design criterion.  So far, no general method to design
  good-performing constellation with large diversity for any
  number of transmit antennas and any transmission rate exists.

  In this paper we introduce numerical methods which allows one
  to construct codes with near optimal diversity sum. We also
  demonstrate how these numerical techniques lead to codes with
  excellent diversity product and more generally diversity
  function. The numerical design methods work for any dimensional
  constellation having any transmission rate.
  Codes can be optimized for any signal to noise ratio.
\end{abstract}

\newpage
\Section{Introduction and Model}

One way to acquire reliable transmission with high transmission
rate on a wireless channel is to use multiple transmit or receive
antennas. Either because of rapid changes in the channel
parameters or because of limited system resources, it is
reasonable to assume that both the transmitter and the receiver
don't know about the channel state information (CSI), i.e. the
channel is non-coherent.

In~\cite{ho00a},  Hochwald  and Marzetta study unitary space-time
modulation. Consider a wireless communication system with $M$
transmit antennas and $N$ receive antennas operating in a Rayleigh
flat-fading channel. We assume time is discrete and at each time
slot, signals are transmitted simultaneously from the $M$
transmitter antennas. We can further assume that the wireless
channel is quasi-static over a time block of length $T$.

A signal constellation $\V:=\{ \Phi_1,\ldots, \Phi_L\}$ consists
of $L$ matrices having size $T \times M$ and satisfying $T \ge M$
and $\Phi_k^* \Phi_k = I_M$. The last equation simply states that
the columns of $\Phi_k$ form a ``unitary frame'', i.e. the column
vectors all have unit length in the complex vector space
$\mathbb{C}^T$ and the vectors are pairwise orthogonal.  The
scaled matrices $\sqrt{T} \Phi_k$, $k=1,2,\cdots,L$, represent the
code words used during the transmission. One can verify that the
transmission rate is determined by $L$ and $T$:

$$\mathtt{R}=\frac{\log_2(L)}{T}.$$

Denote by $\rho$ the signal to noise ratio (SNR).  The basic
equation between the received signal $R$ and the transmitted
signal $\sqrt{T} \Phi$ is given through:
$$
R=\sqrt{\frac{\rho T}{M}}\Phi H+W,
$$
where the $M \times N$ matrix $H$ accounts for the multiplicative
complex Gaussian fading coefficients and the $T \times N$ matrix
$W$ accounts for the additive white Gaussian noise. The entries
$h_{m,n}$ of the matrix $H$ as well as the entries $w_{t,n}$ of
the matrix $W$ are assumed to have a statistically independent
normal distribution $\mathcal{CN}(0,1)$. In particular it is
assumed that the receiver does not know the exact values of either
the entries of $H$ or $W$ (other than their statistical
distribution).

The decoding task asks for the computation of the most likely
sent code word $\Phi$ given the received signal $R$.
Denote by $||\ \ ||_F$ the Frobenius norm of a matrix. If
$A=(a_{i,j})$ then the Frobenius norm is defined through  $|| A
||_F=\sqrt{\sum_{i,j} |a_{i,j}|^2}.$ Under the assumption of
above model the maximum likelihood (ML) decoder will have to
compute:
$$
\Phi_{ML}=\displaystyle \arg \max_{\Phi_l \in
\{\Phi_1,\Phi_2,\cdots,\Phi_L\}} {\|R^*\Phi_l\|}_F
$$
for each received signal $R$. (See~\cite{ho00a}).

Let $\delta_m(\Phi_l^* \Phi_{l'})$ be the $m$-th singular value
of $\Phi_l^* \Phi_{l'}$. It has  been shown in~\cite{ho00a} that the pairwise
probability of mistaking $\Phi_l$ for $\Phi_{l'}$ using Maximum
Likelihood Decoding satisfies:

\begin{eqnarray}
P_{\Phi_l,\Phi_{l'}} &=& \mbox{Prob}\left(\mbox{ choose }\Phi_{l'}\mid
\Phi_{l}\mbox{ transmitted } \right)(\rho)\nonumber\\
 &=& \mbox{Prob}\left(\mbox{ choose }\Phi_{l}\mid
\Phi_{l'}\mbox{ transmitted } \right)(\rho) \nonumber\\
&=& \frac{1}{4\pi} \int_{-\infty}^{\infty}\frac{4}{4w^2+1}\prod_{m=1}^M
\left[1+
\frac{(\rho T/M)^2 (1-\delta_m^2 (\Phi_l^*
\Phi_{l'})) }{4(1+\rho T/M)} (4w^2+1)\right]^{-N}\!\! dw \label{exactf}\\
 &\le& \frac{1}{2} \prod_{m=1}^M
\left[1+
\frac{(\rho T/M)^2(1-\delta_m^2 (\Phi_l^*
\Phi_{l'}))}{4(1+\rho T/M)} \right]^{-N}. \label{mainf}
\end{eqnarray}

It is a basic design objective to construct constellations $\V=\{
\Phi_1,\ldots, \Phi_L\}$ such that the pairwise probabilities
$P_{\Phi_l,\Phi_{l'}}$ are as small as possible. Mathematically we
are dealing with an optimization problem with unitary constraints:

\vspace{0.3cm}

Minimize $\displaystyle \max_{l \neq l'} P_{\Phi_l,\Phi_{l'}}$
with the constraints $\Phi_i^*\Phi_i=I$ where $i=1,2,\cdots,L$.

\vspace{0.3cm}

Formula\eqr{mainf} is sometimes referred to as ``Chernoff's
bound''.  This formula is easy to work with, the exact
formula\eqr{exactf} is in general not easy to work with, although
as we demonstrate in this paper it is very useful in the
numerical search of good constellations as well.  Researchers
have been searching for constructions where the maximal pairwise
probability of $P_{\Phi_l,\Phi_{l'}}$ is as small as possible. Of
course the pairwise probabilities depend on the chosen signal to
noise ratio $\rho$ and the construction of constellations has
therefore to be optimized for particular values of the SNR.

The design objective is slightly simplified if one assumes that
transmission operates at high signal to noise ratios.
In~\cite{ho00}, a design criterion for high SNR  is presented and
the problem has been converted to the design of a finite set of
unitary matrices whose pairwise diversity product is as large as
possible. In this special situation several
researchers~\cite{al98,ta00a,sh01,sh02}   came up
with algebraic constructions and we will say more about this in
the next section.

The main purpose of this paper is to develop numerical procedures
which allow one to construct unitary constellations with excellent
diversity for  any set of parameters $M,N,T,L$ and for any signal
to noise ratio $\rho$.

The paper is structured as follows. In
Section~\ref{Sec-diversity} we introduce the diversity function
of a constellation. This function depends on the signal to noise
ratio and it gives for each value $\rho$ an indication how good
the constellation $\V$ will perform. For large values of $\rho$
the diversity function is governed by the diversity product, for
small values of $\rho$ it is governed by the diversity sum. These
concepts are introduced in Section~\ref{Sec-diversity} as
well. The introduced concepts are illustrated on some well known
constellations previously studied in the literature.

Section~\ref{Sec-numeric} is probably the main section of this
paper. In this section we demonstrate how one can use numerical
tools to derive excellent constellations for any set of
parameters $M,N,T,L$ and $\rho$. The main numerical techniques we
are using are the {\em Simulated Annealing Algorithm} and the
{\em Genetic Algorithm} which are both explained in
Section~\ref{Sec-numeric}.

In Subsection~\ref{Sec-2dim} we present a larger set of
2-dimensional code designs and we indicate their performances
through the computation of their diversity and through
simulations.  Subsection~\ref{Sec-ndim} gives results for general
$n$-dimensional constellations. The simulations indicate that in
the design of codes more attention should be given to the
diversity sum (more generally diversity function) which previously
has not been studied.

\Section{The Diversity Function, the Diversity Sum and the
  Diversity Product}                         \label{Sec-diversity}

In this paper we will be concerned with the numerical construction
of constellations where the right hand sides in\eqr{exactf}
and\eqr{mainf}, maximized over all pairs $l,l'$ is as small as
possible for fixed numbers of $T,M,N,L$. As already mentioned this
tasks depends on the signal to noise ratio the system is
operating.  For this purpose we define the {\em exact
  diversity function} dependent on the constellation $\V=\{
\Phi_1,\ldots, \Phi_L\}$ and a particular SNR $\rho$ through:
\begin{equation}                  \label{exact-div}
\mathcal{D}_e(\V,\rho):=
\max_{l \ne l'} \mbox{Prob}\left(\mbox{ choose }\Phi_{l'}\mid
\Phi_{l}\mbox{ transmitted } \right)(\rho)
\end{equation}

For a particular constellation with a large number $L$ of
elements, with many transmit and receive antennas the function
$\mathcal{D}_e(\V,\rho)$ is very difficult to compute. Indeed for
each pair $\Phi_{l'},\Phi_{l}$ it is required to compute the
singular values of the $M\times M$ matrix $\Phi_l^* \Phi_{l'}$ and
then one has to evaluate up to $L(L-1)/2$ integrals of the
form\eqr{exactf} and this has to be done for each value of $\rho$.
Although this task is formidable it can be done in cases where
$T,M,L$ are all in the single digits using e.g. Maple.

Using Chernoff's bound\eqr{mainf} we define a simplified function
called the {\em diversity function} through:
\begin{equation}                  \label{div}
\mathcal{D}(\V,\rho):=
\max_{l \ne l'}
\frac{1}{2} \prod_{m=1}^M
\left[1+
\frac{(\rho T/M)^2}{4(1+\rho T/M)} (1-\delta_m^2 (\Phi_l^*
\Phi_{l'}))\right]^{-N}.
\end{equation}

The computation of $\mathcal{D}(\V,\rho)$ does not require the
evaluation of an integral and the computation requires essentially the
computation of $ML(L-1)/2$ singular values.

\begin{figure}[ht]
\centerline{\psfig{figure=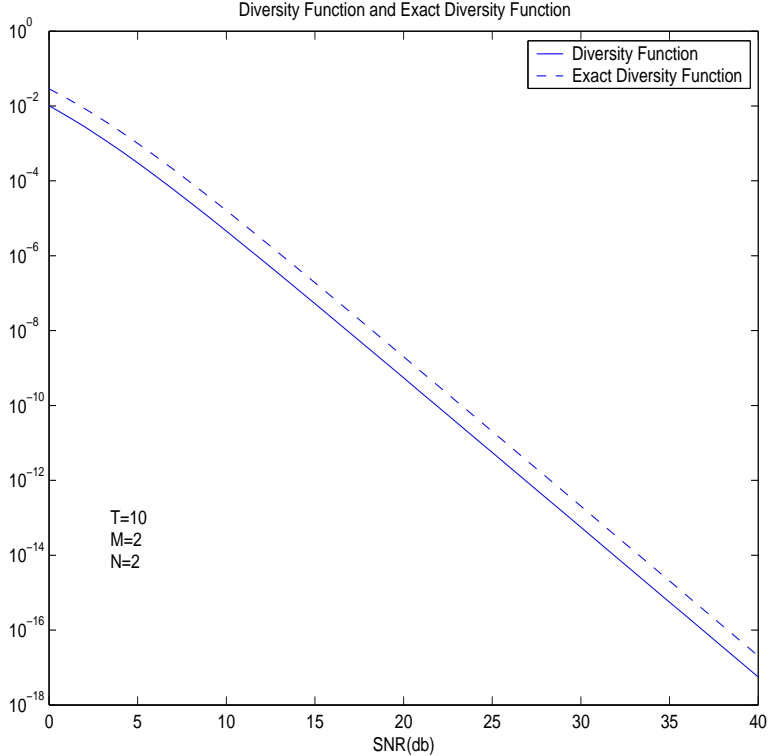,width=4in,height=4in}}
\caption{Diversity function $\mathcal{D}(\V,\rho)$ and exact
  diversity function $\mathcal{D}_e(\V,\rho)$ of a fully
  isotropic constellation.} \label{fig-2}
\end{figure}

The singular values $\delta_m (\Phi_l^*\Phi_{l'})$ are by
definition all real numbers in the interval $[0,1]$ as we assume
that the columns of $\Phi_l,\Phi_{l'}$ form both orthonormal
frames.  The functions $\mathcal{D}_e(\V,\rho)$ and
$\mathcal{D}(\V,\rho)$ are smallest if the singular values
$\delta_m (\Phi_l^*\Phi_{l'})$ are as small as possible. These
numbers are all equal to zero if and only if the column spaces of
$\Phi_l,\Phi_{l'}$ are pairwise perpendicular. We call such a
constellation {\em fully isotropic}.  Since the columns of
$\Phi_l$ generate a $M$-dimensional subspace this can only happen
if $L\leq T/M$. On the other hand if $L\leq T/M$ it is easy to
construct a constellation where the singular values of
$(\Phi_l^*\Phi_{l'})$ are all zero. Just pick $LM$ different
columns from a $T\times T$ unitary matrix. Figure~\ref{fig-2}
depicts the functions $\mathcal{D}_e(\V,\rho)$ and
$\mathcal{D}(\V,\rho)$ for a fully isotropic constellation with
$T=10$ and $M=N=2$.

In order to study the function $\mathcal{D}(\V,\rho)$ more
carefully let
\begin{equation}                              \label{tilderho}
\tilde{\rho}:=\frac{(\rho T/M)^2}{4(1+\rho T/M)}.
\end{equation}
In some small interval $[\rho_1,\rho_2]$ the maximum
in\eqr{div} is achieved for some fixed indices $l,l'$ and in
terms of $\tilde{\rho}$ the function $\mathcal{D}(\V,\rho)$ is of
the form:
$$
\mathcal{D}(\V,\tilde{\rho})=\frac{1}{2\left(
    1+c_1\tilde{\rho}+\cdots +c_M\tilde{\rho}^M\right)^N},
$$
where the coefficients $c_1,\ldots,c_M$ depend on the
particular constellation and on the chosen interval
$[\rho_1,\rho_2]$. For an interval close to zero the dominating
term will be the coefficient~$c_1$. Up to some factor this term
will define the {\em diversity sum} of the constellation. When
$\tilde{\rho}>>0$ then the dominating term will be the
coefficient $c_M$ and up to some scaling this term will define
the {\em diversity product} of the constellation. A constellation
will have a small diversity function for small values of $\rho$
(and presumably performs well in this range) when the
constellation is chosen having a large diversity sum.  A
constellation will have a small diversity function for large
values of $\rho$ (and presumably performs well in this range)
when the constellation is chosen having a large diversity
product.  In the next two subsections we will study the limiting
behavior of $\mathcal{D}(\V,\rho)$ as $\rho$ goes to zero and to
infinity.

\subsection{Design criterion for high SNR}

When the SNR $\rho$ is very large then $\mathcal{D}(\V,\rho)$ can
be approximated via:
\begin{equation}
\mathcal{D}(\V,\rho)\simeq
\max_{l \ne l'}
\frac{1}{2}
\left(
\frac{(\rho T/M)^2}{4(1+\rho T/M)}
\right)^{-NM}
\prod_{m=1}^M
\frac{1}{(1-\delta_m^2 (\Phi_l^*
\Phi_{l'}))^{N}}.
\end{equation}

It is the design objective to construct a constellation $\Phi_1,
\Phi_2,\cdots, \Phi_n$ such that
$$
\min_{l \ne l'} \prod_{m=1}^M(1-\delta_m^2(\Phi_l^*
\Phi_{l'}))
$$
is as large as possible.  This last expression defines in
essence the diversity product.  In order to compare different
dimensional constellations it is customary to use the definition:
\begin{de} (See~\cite{ho00})                             \label{div-prod}
  The {\em diversity product} of a unitary constellation $\V$ is
  defined as
  $$
  \prod \V = \min_{l \ne l'} \left(\prod_{m=1}^M (1-
    \delta_m(\Phi_l^* \Phi_{l'}) ^2)\right)^{\frac{1}{2M}}.
  $$
\end{de}

An important special case occurs when $T=2M$. In this situation
it is customary to represent all unitary matrices $\Phi_k$ in the
form:
\begin{equation}                         \label{specialform}
\Phi_k=\frac{\sqrt{2}}{2} \left(\begin{array}{c}
    I\\
    \Psi_k
  \end{array}\right).
\end{equation}
Note that by definition of $\Phi_k$ the matrix $\Psi_k$ is a $M
\times M$ unitary matrix. The diversity product as defined in
Definition~\ref{div-prod} has then a nice form in terms of the
unitary matrices. For this let $\lambda_m$ be the $m$th
eigenvalue of a matrix. Then we have:
$$
1-\delta_m^2(\Phi_{l'}^* \Phi_l)=\frac{1}{4}
\lambda_m(2I_M-\Phi_l^* \Phi_{l'}-\Phi_{l'}^* \Phi_l)
=\frac{1}{4}\delta_m^2(I_M-\Psi_{l'}^*
\Psi_l)=\frac{1}{4}\delta_m^2(\Psi_{l'}-\Psi_l).
$$
So we have
$$
\prod_{m=1}^M (1-\delta_m^2(\Phi_{l'}^*
\Phi_l))^{\frac{1}{2M}}=\frac{1}{2}\prod_{m=1}^M
\delta_m(\Phi_{l'}-\Phi_l)^{\frac{1}{M}}=\frac{1}{2}|\det(\Psi_{l'}-\Psi_l)|^{\frac{1}{M}}.
$$
When $T=2M$ and the constellation $\V$ is defined as above,
then the formula of the diversity product assumes the simple
form:
\begin{equation}  \label{diversity}
\prod \V=\frac{1}{2} \min_{0 \leq l < l' \leq L}
|\det(\Psi_l-\Psi_{l'})|^{\frac{1}{M}}.
\end{equation}

We call a constellation $\V$ a fully diverse constellation if
$\prod \V > 0$.

A lot of efforts have been taken to construct constellations with
large diversity product. (See
e.g.~\cite{ho00,li02,ha01p,ha02p,sh01,sh02,ta00a}). For the
particular situation $T=2M$ the design asks for the construction
of a discrete subset $\V=\{\Psi_1,\ldots,\Psi_L\}$ of the set of
$M\times M$ unitary matrices $U(M)$. When this discrete subset has
the structure of a discrete subgroup of $U(M)$ then the condition
that $\V$ is fully diverse is equivalent with the condition that
the identity matrix is the only element of $\V$ having an
eigenvalue of 1. In other words the constellation $\V$ is required
to operate fixed point free on the vector space $\mathbb{C}^M$.
Using a classical classification result of fixed point free
unitary representations by Zassenhaus~\cite{za36} Shokrollahi et
al.~\cite{sh01,sh02} were able to study the complete list of fully
diverse finite group constellations inside the unitary group
$U(M)$. Some of these constellations have the best known diversity
product for given fixed parameters $M,N,L$. Unfortunately the
possible configurations derived in this way is somehow limited.
The constellations are also optimized for the diversity product
and as we demonstrate in this paper equal attention should be
given to the diversity sum.

\begin{rem}
  In most of the literature mentioned above researchers focus
  their attention to constellations having the special
  form\eqr{specialform}. Unitary differential modulation is
  used to avoid sending the identity (upper part of every
  element in the constellation) redundantly. This increases the
  transmission rate by a factor of 2 to:
  $$
  \mathtt{R}=\frac{\log_2(L)}{M}=2\frac{\log_2(L)}{T}.
  $$
  Because of this reason we will also focus ourselves in the
  later part of the paper to the special form\eqr{specialform} as
  well.  Nonetheless it will become obvious that the numerical
  techniques also work in the general situation.
\end{rem}

\subsection{Design criterion for low SNR channel}

As we mentioned before a constellation with a large diversity sum
will have a small diversity function at small values of the
signal to noise ratio. This is particularly suitable when the
system operates in a very noisy environment. When $\rho$ is
small, using Formula\eqr{tilderho}, one has the following
expansion:
\begin{multline*}
\prod_{m=1}^M [1+ \frac{(\rho T/M)^2}{4(1+\rho T/M)} (1-\delta_m^2
(\Phi_l^* \Phi_{l'}))]= \prod_{m=1}^M [1+ \tilde{\rho}-\tilde{\rho}
\delta_m^2(\Phi_l^* \Phi_{l'})]\\
=(1+\tilde{\rho})^M-\tilde{\rho} \sum_{m=1}^M \delta_m^2(\Phi_l^*
\Phi_{l'})(1+\tilde{\rho})^{M-1}+O(\tilde{\rho}^2).
\end{multline*}

When $\rho \rightarrow 0$, i.e. $\tilde{\rho} \rightarrow 0$, we
can omit the higher order terms $O(\tilde{\rho}^2)$ and the upper
bound of $P_{\Phi_l,\Phi_{l'}}$ requires that
$$
\sum_m \delta_m^2={\|\Phi_l^* \Phi_{l'} \|}_F^2
$$
is small.  In order to lower the pairwise error probability,
it is the objective to make ${\|\Phi_l^* \Phi_{l'} \|}_F^2$ as
small as possible for every pair of $l, l'$. It follows that at
high SNR, the probability primarily depends on $\prod_{m=1}^M
(1-\delta_m^2)$, but at low SNR, the probability primarily
depends on $\sum_{m=1}^M \delta_m^2$.

 In order to be able to
compare the constellation of different dimensions, we define:
\begin{de}                             \label{div-sum}
The {\em diversity sum} of a unitary constellation $\V$ is defined
as
$$
\sum \V = \min_{l \ne l'} \sqrt{1-\frac{
{\|\Phi_l^*\Phi_{l'}\|}_F^2}{M}}.
$$
\end{de}

Again one has the important special case where $T=2M$ and the
matrices $\Phi_k$ take the special form\eqr{specialform}. In this
case one verifies that
\begin{multline*}
{\|\Phi_{l}^*\Phi_{l'}\|}_F^2=\frac{1}{4}{\|I+\Psi_l^*\Psi_{l'}\|}_F^2
=\frac{1}{4}\tr((I+\Psi_{l'}^*\Psi_l)(I+\Psi_l^*\Psi_{l'}))\\
=\frac{1}{4}\tr(2I+\Psi_{l'}^*\Psi_l+\Psi_l^*\Psi_{l'})=
\frac{1}{4}(4M-(2M-\tr(\Psi_{l'}^*\Psi_l+\Psi_l^*\Psi_{l'})))\\
=\frac{1}{4}(4M-\tr((\Psi_l-\Psi_{l'})^*(\Psi_l-\Psi_{l'})))=
\frac{1}{4}(4M-{\|\Psi_l-\Psi_{l'}\|}_F^2)
\end{multline*}

For the form\eqr{specialform} the diversity sum assumes the
following simple form:
\begin{equation}                                         \label{T2M-div-sum}
\sum \V = \min_{l,l'}\frac{1}{2\sqrt{M}}{\|\Psi_l- \Psi_{l'} \|_F}.
\end{equation}

\begin{rem}
  Without mentioning the term the concept of diversity sum was
  used in~\cite{ho00a1}.  Liang and Xia~\cite[p. 2295]{li02}
  explicitly defined the diversity sum in the situation when
  $T=2M$ using equation\eqr{T2M-div-sum}.
  Definition~\ref{div-sum} naturally generalizes the definition
  to arbitrary constellations.
\end{rem}

As  the formulas make it clear the diversity sum and the
diversity product are in general very different. There is however
an exception. When $T=4$, $M=2$ and the constellation $\V$ is in
the special\eqr{specialform}. If in addition all the $2\times 2$
matrices $\{ \Psi_1,\ldots, \Psi_L\}$
are a subset of the special unitary group
$$
SU(2)=\{ A\in \C^{2\times 2}\mid A^*A=I\mbox{ and }\det A=1\}
$$
then it turns out that the diversity product $\prod \V$ and
the diversity sum $\sum \V$ of such a constellation are the same.
For this note that elements $\Psi_l,\Psi_{l'}$ of $SU(2)$ have
the special form:
$$
\Psi_l=\vier{a}{b}{-\bar{b}}{\bar{a}},\ \Psi_{l'}=\vier{c}{d}{-\bar{d}}{\bar{c}}.
$$
Through a direct calculation one verifies that
$\det(\Psi_l-\Psi_{l'})=|a-c|^2+|b-d|^2$ and
${\|\Psi_l-\Psi_{l'}\|}_F^2=2(|a-c|^2+|b-d|^2)$. But this means
that $\prod \V= \sum \V$ for constellations inside $SU(2)$.

\subsection{Three illustrative examples}    \label{SubSec-I}

The diversity sum and the diversity product govern the diversity
function at low SNR respectively at high SNR. Codes optimized at
these extreme values of the SNR-axis do not necessarily perform
very well on the ``other side of the spectrum''. In this
subsection we illustrate the introduced concepts on three
examples. All examples have about equal parameters, namely $T=4$,
$M=2$ and the size $L$ is 121 respectively 120. The first two
examples are well studied examples from the literature. The third
example is one we derived by numerical methods.

\paragraph{Orthogonal Design:} This constellation has been
considered by several authors~\cite{al98,sh01}. For our purpose
we simply define this code as a subset of $SU(2)$:
$$
\left\{\frac{\sqrt{2}}{2}\left(\begin{array}{cc}
            e^{\frac{2m\pi i}{11}}&e^{\frac{2n\pi i}{11}}\\
            -e^{-\frac{2n\pi i}{11}}&e^{-\frac{2m\pi i}{11}}
            \end{array} \right)|m,n=0,1,\cdots,10\right\}.
$$
The constellation has 121 elements and the diversity sum and the
diversity product are the same.

\paragraph{Unitary Representation of $SL_2(\mathbb{F}_5)$:}
Shokrollahi e. a.~\cite{sh01} derived a constellation using the
theory of fixed point free representations whose diversity
product is near optimal. This constellation appears as a unitary
representation of the finite group $SL_2(\mathbb{F}_5)$ and we
will refer to this constellation as the
$SL_2(\mathbb{F}_5)$-constellation. The finite group
$SL_2(\mathbb{F}_5)$ has 120 elements and this is also the size
of the constellation. In order to describe the constellation let
$\eta = e^{\frac{2\pi i}{5}}$ and define
$$
P=\frac{1}{\sqrt{5}} \left(\begin{array}{cc}
    \eta^2-\eta^3&\eta^1-\eta^4\\
    \eta^1-\eta^4&\eta^3-\eta^2\\
 \end{array} \right),\ \
Q=\frac{1}{\sqrt{5}} \left(\begin{array}{cc}
    \eta^1-\eta^2&\eta^2-\eta^1\\
    \eta^1-\eta^3&\eta^4-\eta^3\\
\end{array} \right).
$$
Then the constellation is given by the set of matrices
$(PQ)^jX$,where $j = 0,1,\cdots,9$, $X$ runs over the set
\begin{multline*}
\{ I_2, P, Q, QP, QPQ, QPQP, QPQ^2, QPQPQ, QPQPQ^2, \\
QPQPQ^2P, QPQPQ^2PQ, QPQPQ^2PQP \}.
\end{multline*}

\begin{figure}[ht]
\centerline{\psfig{figure=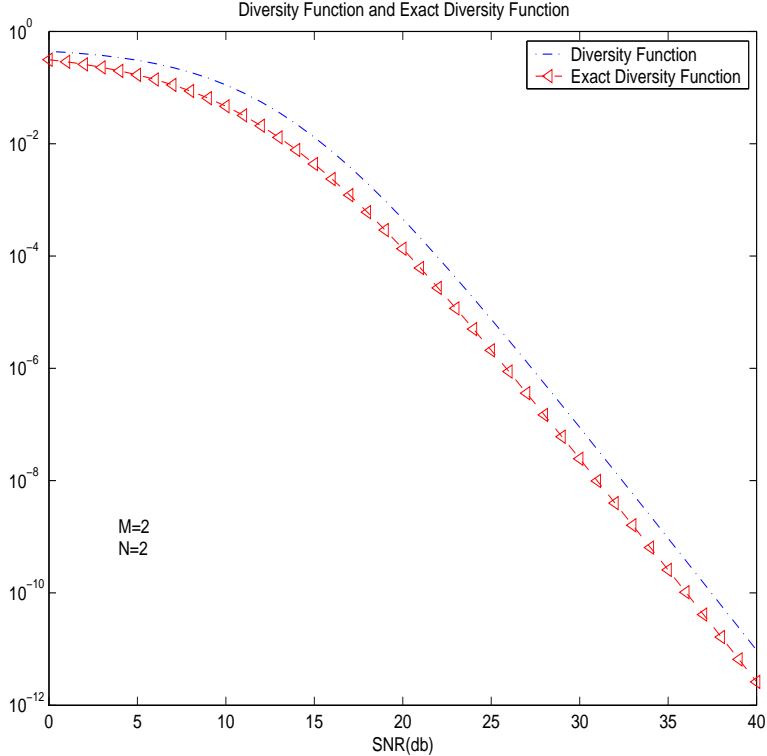,width=4in,height=4in}}
\caption{Diversity function $\mathcal{D}(\V,\rho)$ and exact
  diversity function for the group constellation
$SL_2(\mathbb{F}_5)$.}\label{fig-1}
\end{figure}

The constellation has rate $R = 3.45$ and $\prod
{SL_2(\mathbb{F}_5)}=\sum {SL_2(\mathbb{F}_5)} =
\frac{1}{2}\sqrt{\frac{(3-\sqrt{5})}{2}} \sim 0.309$. The
diversity product of this constellation is truly outstanding. For
illustrative purposes we plotted in Figure~\ref{fig-1}  the
diversity functions of this constellation.

\paragraph{Numerically Derived Constellation:}
Using Simulated Annealing Algorithm we found after short
computation a constellation with near perfect diversity sum. The
constellation is given through a set of 121 matrices
\begin{multline*}
\left\{\Psi_{k,l}:=A^kB^l|A=\left(\begin{array}{cc}
                 -0.9049 + 0.3265*i&   0.1635 + 0.2188*i\\
                 0.0364 + 0.2707*i&  -0.8748 + 0.4002*i
                 \end{array}\right),\right. \\
\left.
B=\left(\begin{array}{cc}
        -0.1596 + 0.9767*i&  -0.1038 + 0.0994*i\\
        0.0833 - 0.1171*i&  -0.9432 + 0.2995*i
        \end{array}\right), k,l=0,1,\cdots,10 \right\}.
\end{multline*}
\bigskip

\begin{ta}
The following table summarizes the parameters of the three
constellations:
\begin{center}
\begin{tabular}{|c|c|c|c|}
  \hline
   & \begin{tabular}{c}
Orthogonal\\design\end{tabular}
 & $SL_2(\mathbb{F}_5)$
& \begin{tabular}{c}
Numerically\\ derived\end{tabular} \\
  \hline
  Number of elements& 121 & 120 & 121\\
  \hline
  diversity sum&  0.1992 &0.309 & 0.3886  \\
  \hline
  diversity product & 0.1992&0.309 &  0.0278 \\ \hline
\end{tabular}
\end{center}
\end{ta}

Of course we were curious about the performance of these three
different codes. Figure~\ref{fig-3} provides simulation results
for each of the three constellations:

\begin{figure}[ht]
\centerline{\psfig{figure=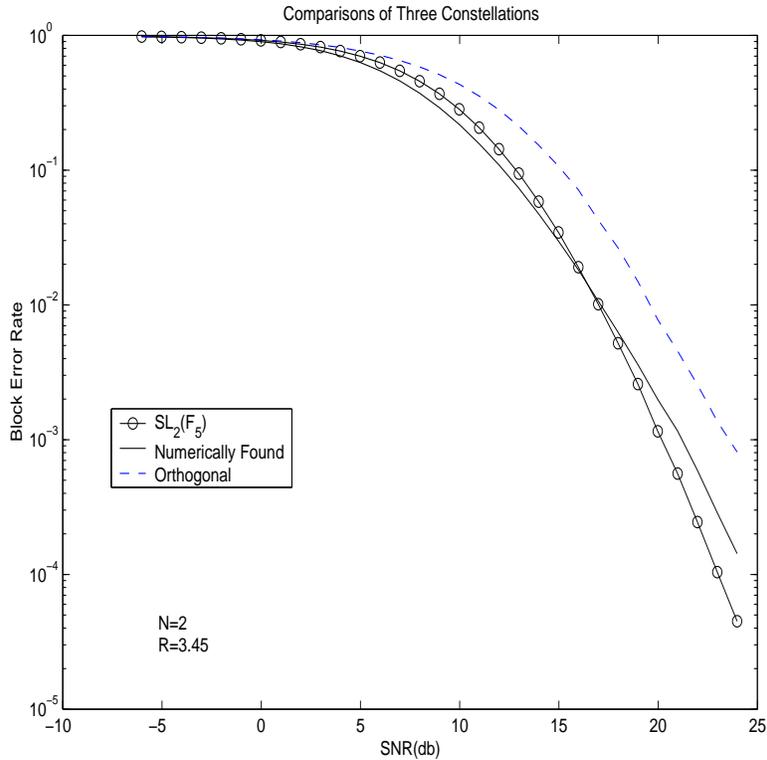,width=4in,height=4in}}
\caption{Simulations of three constellations having sizes
  $T=4$, $M=2$ and $L=120$ respectively $L=121$.}
\label{fig-3}
\end{figure}

Note that the numerically designed code who has a very bad
diversity product is performing very well nevertheless due to the
exceptional diversity sum. These simulation results give an
indication that the diversity sum is a very important parameter
for a constellation.

\Section{Numerical Design of Unitary Constellations with Good
  Diversity}     \label{Sec-numeric}

In order to numerically design constellations with good diversity
it will be necessary to have a good parameterization for the set
of all constellations having size $L$, block length $T$ and
operating with $M$ transmit antennas. In this section we show how
one can use the theory of complex Stiefel manifolds and the
classical Cayley transform to obtain such a parameterization in
all cases.

\subsection{The complex Stiefel manifold}

\begin{de}
  The subset of $T\times M$ complex matrices
  $$
  \S_{T,M}:=\left\{ \Phi\in\C^{T\times M}\mid \Phi^* \Phi =
    I_M\right\}
  $$
  is called the {\em complex Stiefel manifold}.
\end{de}

{}From an abstract point of view a constellation $\V:=\{
\Phi_1,\ldots, \Phi_L\}$ having size $L$, block length $T$ and
operating with $M$ antennas can be viewed as a point in the
complex manifold
$$
\mathcal{M}:=\left(\S_{T,M}\right)^L=
\underbrace{\S_{T,M}\times\cdots \times\S_{T,M}}_{\mbox{$L$
    copies}}.
$$
The search for good constellations $\V$ requires hence the
search for points in $\mathcal{M}$ whose diversity is excellent
in some interval $[\rho_1,\rho_2]$.

Stiefel manifolds have been intensly studied in the mathematics
literature since their introduction by Eduard Stiefel some 50
years ago. A classical paper on complex Stiefel manifolds
is~\cite{at60}, a paper with a point of view towards numerical
algorithms is~\cite{ed99}. The major properties are summarized by
the following theorem:
\begin{thm}                                     \label{Stiefel}
  $\S_{T,M}$ is a smooth, real and compact sub-manifold of
  $\C^{MT}=\R^{2MT}$ of real dimension $2TM-M^2$.
\end{thm}

Some of the stated properties will follow from our further
development.

The following two examples give some special cases.

\begin{exmp}
  $$
  \S_{T,1}=\left\{ x\in\C^T\mid ||x||=\sqrt{\sum_{i=1}^M
      x_i\bar{x}_i}=1 \right\}\subset \R^{2T}
  $$
  is isomorphic to the $2T-1$ dimensional unit sphere
  $S^{2T-1}$.
\end{exmp}

\begin{exmp}
  When $T=M$ then $\S_{T,M}=U(M)$, the group of $M\times M$
  unitary matrices. It is well known that the Lie algebra of
  $U(M)$, i.e. the tangent space at the identity element,
  consists of all $M\times M$ skew-Hermitian matrices. This
  linear vector space has real dimension $M^2$, in particular the
  dimension of $U(M)$ is $M^2$ as well.
\end{exmp}

A direct consequence of Theorem~\ref{Stiefel} is:

\begin{co}
  The manifold $\mathcal{M}$ which parameterizes the set of all
  constellations $\V$ having size $L$, block length $T$ and
  operating with $M$ antennas forms a a real compact manifold of
  dimension $2LTM-LM^2$.
\end{co}

As this corollary makes it clear a full search over the total
parameter space is only possible for very moderate sizes of
$M,L,T$. It is also required to have a good parameterization of
the complex Stiefel manifold $\S_{T,M}$ and we will go after this
task next.

The unitary group is closely related to the complex Stiefel
manifold and the problem of parameterization ultimately boils
down to the parameterization of unitary matrices. For this assume
that $\Phi$ is a $T\times M$ matrix representing an element of
the complex Stiefel manifold $\S_{T,M}$. Using Gramm-Schmidt one
constructs a $T\times (T-M)$ matrix $V$ such that the $T\times T$
matrix $\left[ \Phi\mid V\right]$ is unitary. Define two $T\times
T$ unitary matrices $\left[ \Phi_1\mid V_1\right]$ and $\left[
  \Phi_2\mid V_2\right]$ to be equivalent whenever
$\Phi_1=\Phi_2$. A direct calculation shows that two matrices are
equivalent if and only if there is $(T-M)\times (T-M)$ matrix $Q$
such that:
\begin{equation}                                       \label{Q-matrix}
\left[ \Phi_2\mid V_2\right]=\left[ \Phi_1\mid V_1\right]\vier{I}{0}{0}{Q}.
\end{equation}
Identifying the set of matrices $Q$ appearing in\eqr{Q-matrix}
with the unitary group $U(T-M)$ we get the result:
\begin{lem}                                   \label{Lem-par}
  The complex Stiefel manifold $\S_{T,M}$ is isomorphic to the
  quotient group
  $$
  U(T)/U(T-M).
  $$
\end{lem}
This lemma let us verify the dimension formula for $\S_{T,M}$
stated in Theorem~\ref{Stiefel}:
$$
\dim \S_{T,M}=\dim U(T)-\dim U(T-M)=T^2-(T-M)^2=2TM-M^2.
$$

The  section makes it clear that a good parameterization of the
set of constellations $\V$ requires a good parameterization of the
manifold $\mathcal{M}$ and this in turn requires a good
parameterization of the unitary group $U(M)$.

Once one has a nice parameterization of the unitary group $U(M)$
then Lemma~\ref{Lem-par} provides a way to parameterize the
Stiefel manifold $\S_{T,M}$ as well. Parameterizing $U(T)$ modulo
$U(T-M)$ is however an `over parameterization'. Edelman, Arias and
Smith~\cite{ed99} explained a way on how to describe a local
neighborhood of a (real) Stiefel manifold $\S_{T,M}$. The method
can equally well be applied in the complex case. We do not pursue
this parameterization in this paper and leave this for future
work.

In the remainder of this paper we will concentrate on
constellations having the special form\eqr{specialform}. From a
numerical point of view we require for this a good
parameterization of the unitary group and the next subsection
provides an elegant way to do this.

\subsection{Cayley transformation}

There are several ways to represent a unitary matrix in a very
explicit way. One elegant way makes use of the classical Cayley
transformation. In order that the paper is self contained we
provide a short summary. More details are given
in~\cite[Section~22]{pr94} and~\cite{ha02a}.

\begin{de}
  For a complex $M \times M$ matrix $Y$ which has no eigenvalues at $-1$,
  the Cayley transform of $Y$ is defined to be
  $$
  Y^c= (I + Y)^{-1}(I-Y),
  $$
  where $I$ is the $M \times M$ identity matrix
\end{de}
Note that $(I+Y)$ is nonsingular whenever $Y$ has no eigenvalue at -1.
One immediately verifies that $(Y^c)^c=Y$. This is in analogy to
the fact that the linear fractional transformation
$f(z)=\frac{1-z}{1+z}$ has the property that $f(f(z))=z$.

Recall that a matrix $M$ is skew-Hermitian whenever $A^*=-A$. The
set of $M\times M$ skew-Hermitian matrices forms a linear
subspace of $\C^{M\times M}\cong \R^{2M^2}$ having real dimension
$M^2$. This is the Lie algebra of the unitary group $U(M)$. The
main property of the Cayley transformation is summarized in the
following theorem. (See e.g.~\cite{ha02a,pr94}).
\begin{thm}
  When $A$ is a skew-Hermitian matrix then $(I+A)$ is nonsingular
  and the Cayley transform $V:=A^c$ is a unitary matrix. Vice versa
  when $V$ is a unitary matrix which   has no eigenvalues at $-1$
  then the Cayley transform $V^c$ is skew-Hermitian.
\end{thm}

This theorem allows one to parameterize the open set of $U(M)$
consisting of all unitary matrices whose eigenvalues do not
include $-1$ through the linear vector space of skew-Hermitian
matrices.

The Cayley transformation is very important for the numerical
design of constellations because it makes the local topology of
$U(M)$ clear. One can see that most optimization method require
us to consider the neighborhood of one element in $U(M)$.

\subsection{Simulated Annealing Algorithm (SA)}

In our numerical experiments we have considered several methods.
Because there is a large number of target functions the best known
optimization algorithms such as Newton's Methods~\cite{no99,ed99}
and the Conjugate Gradient Method~\cite{no99,ed99} are difficult
to implement.

Two algorithms of very general nature, the {\em Simulated
  Annealing Algorithm} and the {\em Genetic Algorithm} turned out
to be very practical. In this and the next subsection we describe
these algorithms.

Simulated Annealing is a method which mimics the process of
melted metal getting cooled off. In the annealing process of the
melted metal, first the metal is heated to melt, then the
temperature is getting down gradually. The metal will get to a
minimized energy state if the temperature is lowering slow
enough. For more details about this algorithm, we refer
to~\cite{aa89,la87b,ot89}.

In fact, we would rather call it a general method instead of a
concrete algorithm. Generally speaking, for a given optimization
problem, we always take an initial solution in some certain way.
Then consider the ``neighborhood'' of this solution, we will
accept the solution in the ``neighborhood'' according to some
predefined criterion which might include a probability threshold.

SA is a stochastic method which will find better solution up to
the optimal point in a iterative way. And with Cayley
transformation which is a good representation of any dimensional
unitary matrix, one can see that this method can be applied to
any dimensional constellation design based on the algebraic
structure we proposed as above.

The algorithm we are using also depends on the initial guess. In
$2$ dimensional constellation search, one can use the Brute Force
with coarser grid to find a good guess, then use SA to
iteratively find better and better points. But when it comes to
higher dimension, it is not feasible to use this approach anymore
with limited computing power. One way to find a good initial
guess is to start from an existing constellation (such as a group
constellation).

Our implementation of the algorithm can be summarized in the
following way:

\begin{enumerate}
\item Generate initial generators of the whole constellation.
  For low dimension, one can try to use coarser mesh to find
  better starting points, for high dimension one can randomly
  choose an initial guess.

\item Using the Cayley transformation of the generators generate
  randomly a new constellation in the neighborhood of the old
  constellation where the selection is done using a Gaussian
  distribution with decreasing variances as the algorithm
  progresses.

\item Calculate the diversity function (product, sum) of the
  newly constructed constellation.

\item If the new constellation has better diversity function
  (product, sum), then accept the new constellation. If not,
  reject the new constellation and keep the old constellation (or
  accept it according to Metropolis's criterion~\cite{me53}).

\item Check the stopping criterion, if satisfied, then stop,
  otherwise go to $2$ and continue the iteration.
\end{enumerate}

\subsection{Genetic Algorithm (GA)}

Genetic Algorithm~\cite{ho75b} is a optimization algorithm
proposed by J.H. Holland emulating the  evolutionary process of species.
GA doesn't assume much specific information about the
given problem. First a group of candidates are selected in a
certain way and they are encoded using binary coding in most
cases.  Then consider the offsprings of the candidates (the mating
processing is also defined with respect to the specific problem). In
the same time there might be new random candidates added in. The
whole idea is ``better survive'', which really means some
candidates with higher {\em Fitness Evaluation Function} keep staying
in the group and other ``worse'' ones will get discarded. The
algorithm will stop if some threshold is reached after certain
number of iterations. The interested reader will find further
details in~\cite{co99,ma99b}.

Our adapted algorithm is summarized in the following way:

\begin{enumerate}
\item Generate the initial population randomly with the desired
  size.

\item For every individual in the population, calculate its
  Fitness Evaluation Function. In our situation this is the
  minimum distance from this individual to all other individuals.
  Of course this depends on the criterion we want to optimize
  (diversity function, product, sum).

\item Replace the individuals which have the worst Fitness
  Evaluation Function with the same number of randomly chosen
  individuals. (If necessary, randomly choose certain number of
  individuals for mutation).

\item Calculate the diversity function (product, sum) of the
  given size population, if it improves, accept it, otherwise,
  restore the population to the previous population.

\item Check the stopping criterion, if satisfied, then stop,
  otherwise go to 2 and continue the iteration.
\end{enumerate}

\subsection{Constellations with algebraic structure}

As explained in the beginning of this section the number of
parameters to be optimized explodes for even moderately chosen
numbers $T,M$ and $L$. The full parameterization of a
constellation with size $L=120$, $T=4$ and $M=2$ has 1440 free
parameters and this is out of the range for a naive implementation
of the simulated annealing algorithm. The genetic algorithm is
still feasible in this range.

In this subsection we explain how one can restrict the parameter
space to judiciously chosen subsets and do the optimization
inside this restricted parameter space only.

Consider a general constellation of square unitary matrices:
$$
\V=\{\Psi_1, \Psi_2, \cdots, \Psi_L\}.
$$
In order to calculate the diversity product (or sum), one
needs to do $\frac{L(L-1)}{2}$ calculations: $
|\det(\Psi_i-\Psi_j)|$, for every different pair $i,j$.

If one deals with a group constellation then one needs only to
calculate $L-1$ such determinant calculations and this is one of
the big advantages of group constellations. This is a direct
consequence of
$$
|\det(\Psi_i-\Psi_j)|=|\det(\Psi_i)\det(I-\Psi_i^*\Psi_j)|
=|\det(I-\Psi_i^*\Psi_j)|,
$$
where $\Psi_i^*\Psi_j$ is still in the group.

Group constellations are however very restrictive what the
algebraic structure is concerned. In the following we are going to
present some constellations in between general constellations and
group constellations. Our constellations have some small number of
generators. This will ensure that the total parameter space to be
searched is limited as well. Like for group constellations we can
reduce the number of targets. We start with an example:
\begin{exmp}
Consider the  constellation
$$
\V = \{A^kB^l|A, B \in U(M), k=0, \cdots, p, l=0, \cdots, q
\},
$$
The parameter space for this constellation is $U(M)\times
U(M)$, this is a manifold of dimension~$2M^2$ and the number of
elements is $(p+1)(q+1)$. If one has to compute
$|\det(\Psi_i-\Psi_j)|$ for every distinct pair this would
require $\left(\begin{array}{c} (p+1)(q+1) \\ 2
  \end{array}\right)$. We will show in the following that the
same result can be obtained by doing $2pq+p+q$ determinant
computations. This is in analogy to the situation of group
constellations.

Let $\Psi_i$ and $\Psi_j$ be two distinct elements having the
form $A^{k_1}B^{l_1}$ and $A^{k_2}B^{l_2}$ respectively.

We have now several cases. When
 $k_1 =k_2$, then necessarily  $l_1 \neq l_2$ and the distance is
 computed as
$$|\det(A^{k_1}B^{l_1}-A^{k_2}B^{l_2})|=|\det(I-B^{|l_2-l_1|})|,$$
where $|l_2-l_1|$ is an integer between $1$ and $q$.

If $l_1 =l_2$, then we have $k_1 \neq k_2$ and the distance is computed as
$$|\det(A^{k_1}B^{l_1}-A^{k_2}B^{l_2})|=|\det(I-A^{|k_2-k_1|})|,$$
where $|k_2-k_1|$  is an integer between $1$ and $p$.

If $(k_1 < k_2 \;\; \mbox{and} \;\; l_1 < l_2)$ or $(k_1 > k_2
\;\; \mbox{and} \;\; l_1 > l_2)$, we have

$$|\det(A^{k_1}B^{l_1}-A^{k_2}B^{l_2})|=|\det(I-A^{|k_2-k_1|}B^{|l_2-l_1|}),$$
where $1\leq |k_2-k_1|\leq p$ and $1\leq |l_2-l_1|\leq q$.

Similarly if
 $(k_1 < k_2 \;\; \mbox{and} \;\; l_1 > l_2)$ or $(k_1 > k_2
\;\; \mbox{and} \;\; l_1 < l_2)$ then
$$
|\det(A^{k_1}B^{l_1}-A^{k_2}B^{l_2})|=|\det(A^{|k_2-k_1|}-B^{|l_2-l_1|})|,
$$
with $1\leq |k_2-k_1|\leq p$ and $1\leq |l_2-l_1|\leq p$.

The total number of distances to be computed is in total
equal to $2pq+p+q$.
\end{exmp}

Above idea can be applied to describe in an algebraic way
constellations where the number of distances is always
considerably smaller than the total number of comparisons
involving all possible pairs. The following provides a list of
similar generator sets for constellations.

\begin{enumerate}
\item $\V = \{A^kB^l|A, B \in U(M), k=0, \cdots, p, l=0, \cdots,
  q \}$

\item $\V = \{A^kB^lC^m|A, B, C \in U(M), k=0, \cdots, p, l=0,
  \cdots, q, m=0, \cdots, r \}$

\item $I, A, AB, ABA, ABAB, ABABA, \cdots$

\item $I, A, AB, ABC, ABCA, ABCAB, ABCABC, ABCABCA, \cdots $

\item $I, AB, A^2B^2, A^3B^3, \cdots $
\item $I, ABC, A^2B^2C^2, A^3B^3C^3, \cdots $
\end{enumerate}

The last generator set can also be found in~\cite{li02} to design
codes algebraically . Above construction ideas can be applied to
do ``product constructions''. For instance,
\begin{enumerate}\setcounter{enumi}{6}
\item Let $S_1=\{I, C, C^2, C^3, \cdots \}$ and $S_2=\{I, A, AB,
  ABA, \cdots\}$ and consider the cartesian product constellation
  $$
  S = S_1 \times S_2=\{s_1s_2|s_1 \in S_1, s_2 \in S_2\}.
  $$

\item Similarly, let $S_1=\{I, A, AB, ABA, \cdots\}$ and
  $S_2=\{I, C, CD, DCD, CDCD, \cdots \}$ and consider the
  constellation
   $$
  S = S_1 \times S_2=\{s_1s_2|s_1 \in S_1, s_2 \in S_2\}.
  $$
\end{enumerate}

For a $P,P' \in U(M)$ and any matrix $Q \in C^{M \times M}$, we
have
$$
{\|PQP'\|}_F={\|Q\|}_F.
$$
Like for the diversity product it follows from above that the
number of targets to be checked in order to compute the diversity
sum is reduced in a similar way. In fact, it can be shown that
the method works in the same way if one wants to optimize the
diversity function at a certain SNR.  This is essential the case
because:
\begin{equation} \label{sing}
\delta_m(PQP')=\delta_m(Q),
\end{equation}
for $m=1,2,\cdots,M$.

With constellations as above we  are able to reduce the dimension
of  the   parameter space   and in    the same time    we have  a
considerable  reduction in the number  of  targets to be checked.
This is absolutely essential for the Simulated Annealing Algorithm
which we described before.

\begin{rem}
  Because of\eqr{sing} one can use the algebraic structures we
  presented also for the design of constellations where exact
  diversity function is optimized at a certain SNR. The idea is
  appealing but of course it involves the evaluation of an
  improper integral in every loop and this is computationally
  very expensive. For a constellation of size 3 this is however
  easily possible and we have done this in the sequel.

The first constellation is $3$ element $2$ dimensional
constellation with optimal diversity product and diversity
sum(consult Appendix~\ref{three} for more details):

$$
\left\{I, A, B|A=\left(\begin{array}{cc}
      e^{\frac{2\pi i}{3}}&0\\
      0& e^{\frac{-2\pi i}{3}}\\
              \end{array}\right),
            B=\left(\begin{array}{cc}
                e^{\frac{4\pi i}{3}} & 0\\
                0& e^{\frac{-4\pi i}{3}}\\
\end{array}\right)\right\}
$$

The second constellation is optimized at $5$ db by Simulated
Annealing Algorithm based on exact diversity function:
\begin{multline*}
  \left\{I, A_1, B_1|A_1=\left(\begin{array}{cc}
        -0.4530000-0.7804689*i,&0.2197119-0.3706563*i\\
        -0.1733377-0.3944787*i,&-0.5439448+0.7200449*i\\
              \end{array}\right),\right.   \\
          \left. B_1=\left(\begin{array}{cc}
                -0.4475155+0.7358245*i,& -0.1243078+0.4927877*i\\
                0.1862253+0.4728766*i,&-0.5378253-0.6726454*i\\
              \end{array}\right)\right\} .
\end{multline*}
The simulation results are given in Figure~\ref{fig-4}.\clearpage
\begin{figure}[ht]
  \centerline{\psfig{figure=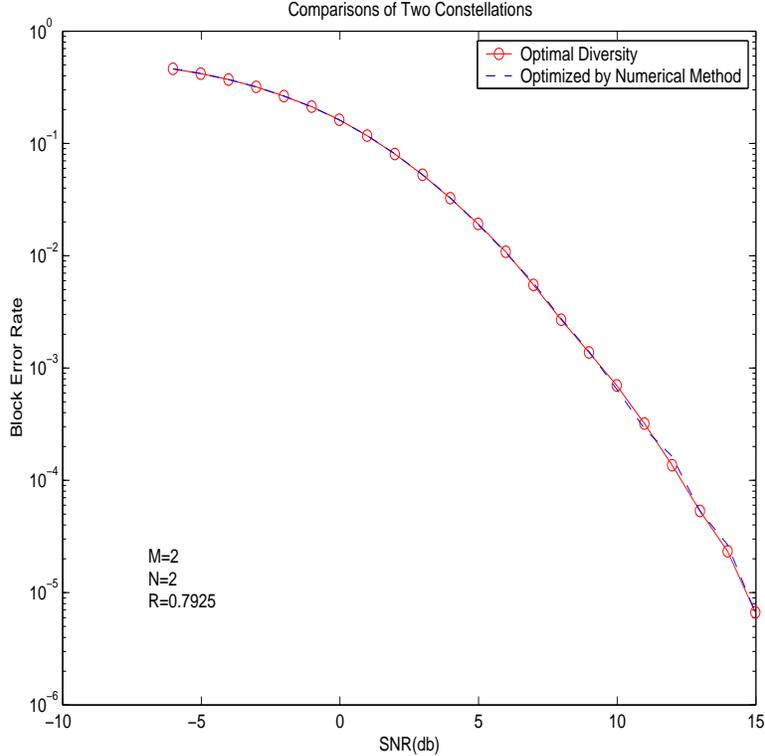,width=4in,height=4in}}
\caption{Optimized constellation by exact diversity function
  versus constellation with optimal diversity product and sum.
  {}From the graph, one can see that the two constellation almost
  have the same performance, except there are minor deviations at
  certain interval. Actually it can be verified that the elements
  in the second constellation are similar to (up to one unitary
  matrix) those in the first constellation, so it has also
  optimal diversity product and diversity sum. For more details
  about $2$ dimensional optimal constellation with $3$ elements,
  one can look at Appendix~\ref{three}. }
\label{fig-4}
\end{figure}
\end{rem}

\Section{Numerical Results of Constellation Design}

\subsection{2 Dimensional constellation design} \label{Sec-2dim}

We are going to talk about $2$ dimensional constellation design
separately because unlike higher dimensional ones, we have a very
explicit representation of $U(2)$:

\begin{equation}  \label{u2}
U(2)=\left\{\left(\begin{array}{cc}
  a & b e^{i\theta} \\
  -\bar{b} & \bar{a} e^{i\theta}
\end{array}\right)||a|^2+|b|^2=1, \qquad 0 \leq \theta \leq 2\pi
\right\} .
\end{equation}

From above, one can see that $U(2)$ is $4$ dimensional, something
we have already mentioned before.  So when we consider the
finitely generated constellation with a good algebraic structure,
it is possible to use Brute Force method to calculate the almost
optimal generator. Of course in this process, finer grid and more
computing power will make better result.

\newpage

\begin{ta}                              \label{Tab-1}
Comparisons of different methods and different parameterizations:
\begin{center}
\begin{tabular}{|c|c|c|c|}
  \hline
    & size&diversity product& structure \\
  \hline
  Simulated Annealing& 125& 0.2127& $A^kB^lC^m$ \\
  \hline
  Simulated Annealing & 120&0.2202& $A^kB^l$\\
  \hline
  Simulated Annealing&121& 0.2417 & $A^kB^l$ \\
  \hline
  Brute Force&120& 0.1914& $A^kB^l$\\
  \hline
   Genetic Algorithm& 120&0.2377 &N/A\\
  \hline
\end{tabular}
\end{center}
\end{ta}
\vspace{0.3cm}

\begin{ta}
  Optimization of the the diversity product using Simulated
  Annealing and different generator sets. All computations were
  accomplished in less than $3$ minutes.
\begin{center}
\begin{tabular}{|c|c|c|c|c|c|c|c|}
\hline
structure/size&    36&       49&       64&            256&400& 900&       10000\\  \hline
$A^kB^l$&            0.3860&   0.3781&   0.2742&        0.1025&
0.0866&    0.0834&    0.0158\\ \hline
$AB$&                0.3205&   0.2659&   0.2450&        0.1030&
0.0800&    0.0579&    0.0122\\ \hline
$A^kB^k$&            0.3769&   0.3502&   0.3090&        0.1651&
0.1342&    0.0820&    0.0187\\ \hline
\end{tabular}

\vspace{0.3cm}

\begin{tabular}{|c|c|c|c|c|c|c|c|}
\hline
structure/size&    27&       64&       216&           343&
512&       729&       9261 \\ \hline
$A^kB^lC^m$&         0.3418&   0.2616&   0.1833&        0.1401&
0.0632&    0.1012&    0.0031\\ \hline
$ABC$&               0.3299&   0.1832&   0.1033&        0.0725&
0.0555&    0.0430&    N/A\\ \hline
$A^kB^kC^k$&         0.4122&   0.2512&   0.0583&        0.0206&
0.0087&    0.0039&    N/A\\ \hline
\end{tabular}
\end{center}
\end{ta}

The above tables show some of the numerical results in $2$
dimension. The constellations in Table~\ref{Tab-1} all have about
120 elements and they are therefore comparable with the
$SL_2(\mathbb{F}_5)$ constellation described in
Subsection~\ref{SubSec-I}. The obtained diversity products are
reasonable good but fall short from the outstanding diversity
product of $0.309$ which the $SL_2(\mathbb{F}_5)$ constellation
has.

The numerical design method seems to be particular powerful if
one seeks constellations in higher dimension or if one wants to
optimize the diversity sum. One can see the performance of higher
dimensional constellations in Subsection~\ref{high}.

\begin{ta} Comparison of different methods or
different parameters
\begin{center}
\begin{tabular}{|c|c|c|c|}
  \hline
    &size& diversity sum &structure\\
  \hline
  Simulated Annealing&125& 0.3919&$A^kB^lC^m$ \\
  \hline
  Simulated Annealing& 120&0.3696&$A^kB^l$ \\
  \hline
  Simulated Annealing&121& 0.3886&$A^kB^l$ \\
  \hline
  Brute Force& 120&0.3673& $A^kB^l$\\
  \hline
   Genetic Algorithm&120& 0.3867&N/A \\
  \hline
\end{tabular}
\end{center}
\end{ta}

\newpage
\begin{ta}
Optimization of the the diversity sum using Simulated
  Annealing and different generator sets. All computations were
  accomplished in less than $3$ minutes.
\begin{center}

\begin{tabular}{|c|c|c|c|c|c|c|c|}
\hline
structure/size&           36&       49&       64&            256&
400&       900&       10000\\ \hline
$A^kB^l$&                   0.5113&   0.4733&   0.4474&
0.2875&     0.2504&    0.1848&    0.0785\\ \hline
$AB$&                       0.5530&   0.4240&   0.3821&
0.1994&     0.1629&    0.1064&    0.0310\\ \hline
$A^kB^k$&                   0.5466&   0.5121&   0.4735&
0.3088&     0.2637&    0.2047&    0.0869\\ \hline
\end{tabular}

\vspace{0.3cm}

\begin{tabular}{|c|c|c|c|c|c|c|c|}
\hline
structure/size&           27&       64&       216&           343&
512&       729&       9261\\ \hline
$A^kB^lC^m$&                0.5400&   0.4210&   0.2992&
0.2663&     0.2099&    0.2060&    0.0772\\ \hline
$ABC$&                      0.5382&   0.4497&   0.2614&
0.2065&     0.1695&    0.1447&    0.0398\\ \hline
$A^kB^kC^k$&                0.5630&   0.4271&   0.2864&
0.2198&     0.1969&    0.1423&    N/A\\ \hline
\end{tabular}
\end{center}
\end{ta}

Numerical method is extremely good for designing constellation
with large diversity sum, most of the results are the best or very
close to the best codes. It can be verified that the
$SL_2(\mathbb{F}_5)$ group code also has diversity sum $0.309$, we
believe that would be the reason its performance is so good. Again
as for higher dimension and large cardinality, numerical methods
are even more powerful.

We have tried SA on constellation with $10000$ elements using the
structure $A^kB^l$.  For this we can get a diversity sum of
$0.1000$. For the construction in~\cite{ha02p} which can be proven
to be asymptotically the best as a subset of $SU(2)$, we can only
get a diversity sum of $0.0654$ with $8433$ elements or a
diversity sum of $0.0604$ with $10770$ elements.
\subsection{Constellation design for any dimension}
\label{Sec-ndim}

As we mentioned before, we can always choose an existing
constellation as staring point as our numerical method. In the
sequel, we use the group constellation $G_{21,4}$ in~\cite{sh01}:

$$
\V_1=\{A^kB^l| A=\left(\begin{array}{ccc}
         \eta&0&0\\
         0&\eta^4&0\\
         0&0&\eta^{16}\\
     \end{array}\right), B=\left(\begin{array}{ccc}
                                  0&1&0\\
                                  0&0&1\\
                                  \eta^7&0&0\\
                                  \end{array}\right),
                                  k=0,1,\cdots,20, l=0,1,2\}
$$

One can verify that

$$\prod \V_1=0.3851$$

It seems like the $G_{21,4}$ is already very good constellation,
our algorithm only improves a little. However one can check most
of the case, the algorithm will improve much compared to the
original group constellation.

$$\V_2=\{A^kB^l|k=0,1,\cdots,20,l=0,1,2 \},$$
where
$$A=\left(\begin{array}{ccc}
         0.9415 + 0.3155*i&0.0573 - 0.0222*i&0.0496 +0.0882*i\\
         0.0160 - 0.0555*i&0.4005 + 0.9136*i&0.0326 - 0.0212*i\\
         0.0579 + 0.0855*i&-0.0312 - 0.0099*i&0.1384 - 0.9844*i\\
                 \end{array}\right),$$
$$ B=\left(\begin{array}{ccc}
        0.0175 + 0.0095*i&0.9997 + 0.0111*i&0.0079 + 0.0042*i\\
        0.0086 + 0.0100*i&-0.0082 + 0.0040*i&0.9999 + 0.0036*i\\
        -0.4836 + 0.8750*i&0.0004 - 0.0198*i&-0.0045 - 0.0126*i\\
        \end{array}\right).$$

One verifies that

$$\prod \V_2=0.3874.$$

Another approach is to choose the starting points randomly. We
have the following tables to show the results of some experiments.
The experiments are based on the computation on a Intel Pentium
800MHz PC and no computation lasts longer than 3 minutes.

\vspace{0.3cm}
\begin{ta}
Numerical results on diversity product using Simulated Annealing
Algorithm

\begin{center}
\begin{tabular}{|c|c|c|c|c|}
  \hline
   \begin{tabular}{c}
   number of\vspace*{-2mm}\\ elements\end{tabular}& dim=2 & dim=3& dim=4& dim=5\\
  \hline
  4 & 0.7071 &0.7657&0.7388&0.6768 \\
  \hline
  9 & 0.5701  &0.5754&0.4774&0.4259\\
  \hline
  16 & 0.4018  &0.4574&0.4651&0.3877\\
  \hline
  25 & 0.3443 &0.3834&0.3809&0.3467\\
  \hline
  36 & 0.2865 &0.3450&0.3501&0.3760 \\
  \hline
\end{tabular}
\end{center}
\end{ta}

\begin{ta}
Numerical results on diversity sum using Simulated Annealing
Algorithm
\begin{center}
\begin{tabular}{|c|c|c|c|c|}
  \hline
   \begin{tabular}{c}
   number of\vspace*{-2mm}\\ elements\end{tabular}& dim=2 & dim=3& dim=4& dim=5\\
  \hline
  4 & 0.8147 &0.8160&0.7861&0.7377 \\
  \hline
  9 & 0.6956  &0.6861&0.6539&0.6389\\
  \hline
  16 & 0.5908 &0.6459&0.6288&0.5916\\
  \hline
  25 & 0.5618 &0.6268&0.6190&0.5795 \\
  \hline
  36 & 0.5286 &0.6054&0.6148&0.5853 \\
  \hline
\end{tabular}
\end{center}

\end{ta}

\vspace{0.3cm}

From above table, SA works well when the cardinality of the
constellation is small. When it comes to large number of
constellation, SA works even better as we mentioned above.  Also
SA's results on higher dimension showing the pretty much the same
size as $2$ dimension which makes us believe it will also work in
higher dimension, although we don't have many comparisons
available.

We have the following tables to show the results of some
experiments using {\em Genetic Algorithm}. The experiments are
based on the computation on a Intel Pentium 800MHz PC and no
computation lasts longer than 3 minutes.

\vspace{0.3cm}

\begin{ta}
Numerical results on diversity product using Genetic Algorithm

\begin{center}
\begin{tabular}{|c|c|c|c|c|c|}
  \hline
   \begin{tabular}{c}
   number of\vspace*{-2mm}\\ elements\end{tabular}& dim=2&optimal DP&dim=3& dim=4& dim=5\\
  \hline
  3 & 0.8644&$\frac{\sqrt{3}}{2}$=0.8860 &0.8264&0.7305&0.6737\\
  \hline
  4 & 0.8051&$\sqrt{\frac{2}{3}}$=0.8165 &0.7343&0.6521&0.6305\\
  \hline
  6 & 0.6924&unknown  &0.6632&0.6154&0.5721\\
  \hline
  10 & 0.5768&unknown &0.5497&0.5742&0.4942\\
  \hline
\end{tabular}
\end{center}
\end{ta}

\newpage

\begin{ta}

Numerical results on diversity sum using Genetic Algorithm

\begin{center}
\begin{tabular}{|c|c|c|c|c|c|c|c|c|}
  \hline
   \begin{tabular}{c}
   number of\vspace*{-2mm}\\ elements\end{tabular}& dim=2&optimal DS & dim=3& dim=4& dim=5\\
  \hline
  3 & 0.8601&$\frac{\sqrt{3}}{2}$=0.8860 &0.8331&0.8118&0.7798\\
  \hline
  4 & 0.8029&$\sqrt{\frac{2}{3}}$=0.8165  &0.7802&0.7757&0.7492\\
  \hline
  6 & 0.7443&$\sqrt{\frac{3}{5}}$=0.7746  &0.7502&0.7293&0.7176\\
  \hline
  10 & 0.6826&$\frac{\sqrt{2}}{2}$=(0.7071) &0.6981&0.6920&0.6817\\
  \hline
\end{tabular}
\end{center}

\end{ta}

Genetic algorithm works extremely well when it comes to small
number constellation. For instance, look at  $3$ elements $2$
dimension diversity product case, one can check the optimal
diversity product is $\frac{\sqrt{3}}{2} \sim 0.8660 $(see
Appendix~\ref{three}), while our results shows $0.8601$. The
results make us believe for higher dimension, maybe we are
approaching the optimal diversity product(sum) too (In the above
table for $3$ or higher dimension constellation, we don't know the
optimal diversity product). One can see other comparisons with
known optimal diversity product(sum)~\cite{li02}. Whereas when it
comes to huge size constellation, our current GA doesn't seems to
work as well as small number constellation because of drastically
increased complexity. Maybe some more complicated mechanism is
needed to add in to handle this problem more subtly.

\subsection{Performance of different dimension constellation}
\label{high}

Basically our numerical method can be applied to design
constellation for every dimension and any transmission rate. The
following graph shows the comparison of three constellations with
different dimensions with $2$ receiver antennas. The first one is
a $2$ dimensional constellation with $3$ elements ($R=0.7925$)
and optimal diversity product $0.8660$ and optimal diversity sum
$0.8660$.  The second constellation is $3$ dimensional
constellation which has $5$ elements ($R=0.7740$) with diversity
product $0.7183$ and diversity sum $0.7454$. The third
constellation is a $4$ dimensional one consisting of $9$ elements
($R=0.7925$) with diversity product 0.5904 and diversity sum
0.6403. Here we used Genetic Algorithm to optimize the diversity
product and the diversity sum at the same time to acquire the
last two constellations.

One can see that around $5$ db, the second constellation surpasses
the first one and is getting better and better as SNR becomes
larger. This can be easily understood that the diversity function
of the first constellation is approximately dominated by
$1/{\rho^4}$ at high SNR, while the diversity function of the
second constellation is  dominated by $1/{\rho^6}$. The same
explanation can be applied to the third constellation's
performance. One can even foresee that higher dimensional
constellations will perform even better and the BER curve will be
sharper than low dimensional ones. It is believable that higher
dimensional constellation ones will have much more diversity gain
compared to lower dimensional ones as for other transmission rate.

\begin{figure}[ht]
\centerline{\psfig{figure=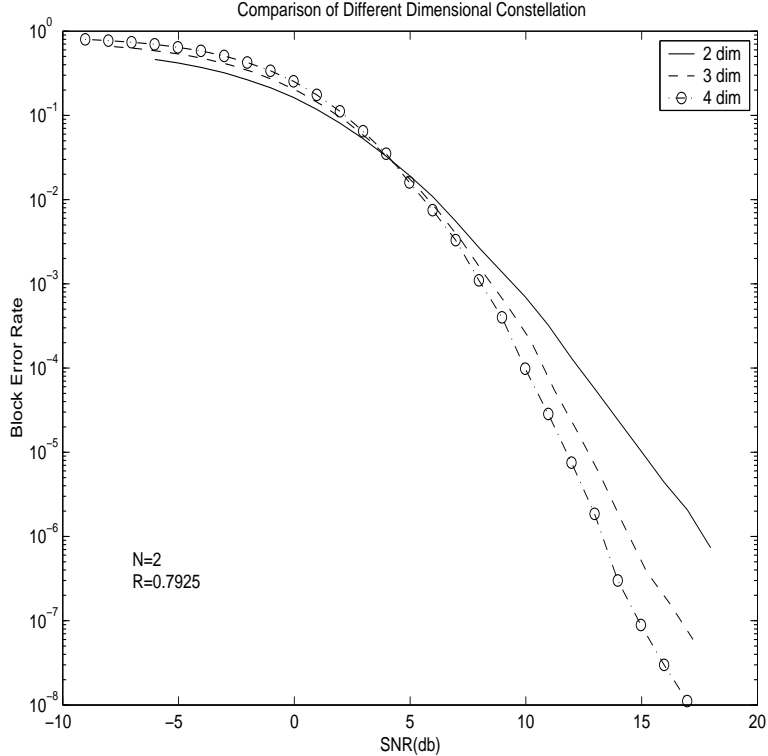,width=4in,height=4in}}
\caption{Performance of different dimensional constellations with
the same rate} \label{fig-5}
\end{figure}

\subsection{Constellations optimized for certain values of SNR by
  diversity function}

Different industrial applications require different level of
reliability of the communication channels. One may want to
optimize the constellation at certain Block Error Rate (BER) or
Signal Noise Ratio (SNR). Algebraically designing codes for this
purpose seems to be too complicated, while the numerical method
can be easily applied for this purpose. In fact as we mentioned
before, algebraic structure is still beneficial for optimizing the
diversity function using simulated annealing algorithm, we also
find our genetic algorithm only needs minor modification to fit in
for this purpose. We have the following graph to show some further
constellations.

\begin{figure}[ht]
\centerline{\psfig{figure=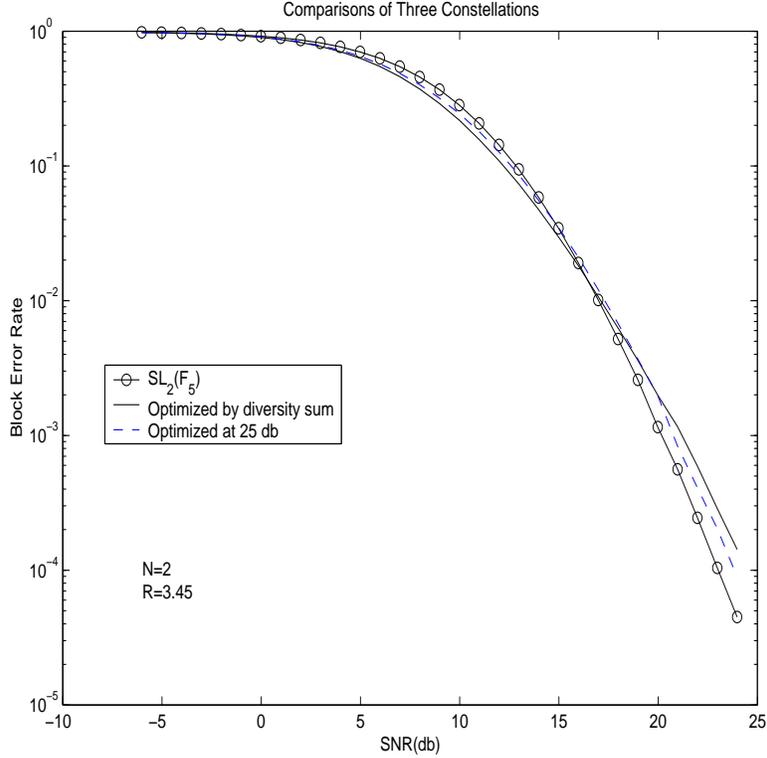,width=4in,height=4in}}
\caption{Optimize the constellations at certain SNR by diversity
function} \label{fig-6}
\end{figure}

We have already seen the first two constellations, the third
constellation is specifically optimized at $25$ db using Simulated
Annealing Algorithm based on structure $A^kB^l$ , one can see that
optimizing at different curve does change the shape of curve.
Consequently one can see the flexibility of the numerical method
employed here.

The third constellation has the following representation:

\begin{multline*}
\left\{A^kB^l|A=\left(\begin{array}{cc}
                 -0.3018715+0.8863567*i&  0.1337423-0.3245896*i\\
                -0.3487261-4.0441897E-02*i& -0.9215508-0.1658271*i
                 \end{array}\right),\right.\\
\left. B=\left(\begin{array}{cc}
        -0.9599144-0.1577551*i&  0.2074585+0.1031435*i\\
        -0.2074221+0.1032166*i&  -0.9598588+0.1580936*i
        \end{array}\right), k,l=0,1,\cdots,10 \right\}.
\end{multline*}

\subsection{General form constellation numerical design}

As above, we discussed numerical methods for the special form of
unitary constellation:

$$\V:=\{ \Phi_1,\ldots, \Phi_L\},$$

where

$$\Phi_k=\frac{\sqrt{2}}{2} \left(\begin{array}{c}
    I\\
    \Psi_k
  \end{array}\right).$$

There are numerical methods available to design general
constellation as well. Right now we are only taking the simplest
approach. In this case, we assume that $T=2M$ and we consider the
following constellation:
$$
\{A^kB|A \in U(M), B=\left(\begin{array}{c}
                  I\\
                  0\\
       \end{array}\right), k=0,1, \cdots, l\}.
$$

With this algebraic structure, whenever we want to calculate the
diversity product or sum or more generally the diversity function, one
needs to calculate $\Phi_{l'}^*\Phi_l$.  It is easily verified that
they will end up being one of the following:

$$
B^*A^kB \qquad  k=-l, \cdots, l.
$$
In this way, we have a finitely generated constellation with
reduced number of targets.

\begin{ta} The following tables show our
  elementary results using SA.

\begin{center}
\begin{tabular}{|c|c|c|c|c|c|c|c|}
  \hline
  size & 3 & 4 & 5 & 6 & 7 & 8 & 9 \\
  \hline
  rate & 0.3962& 0.5000 & 0.5905 & 0.6462& 0.7018 & 0.7500 & 0.7925 \\
  \hline
  diversity sum & 0.8654 & 0.7901 & 0.7889 & 0.7652& 0.7514 & 0.7422 & 0.7369 \\
  \hline
\end{tabular}

\vspace{0.3cm}

\begin{tabular}{|c|c|c|c|c|c|c|c|}
  \hline
  size & 3 & 4 & 5 & 6 & 7 & 8 & 9 \\
  \hline
  rate & 0.3962& 0.5000 & 0.5905 & 0.6462& 0.7018 & 0.7500 & 0.7925 \\
  \hline
  diversity product & 0.8582 & 0.7424 & 0.7330 & 0.6450 & 0.6361 & 0.6216 & 0.5822 \\
  \hline
\end{tabular}
\end{center}
\end{ta}

\Section{Conclusion and Future Work}

In this paper, we study the limiting behavior of the {\em diversity
  function} by either letting the SNR go to infinity or to zero.
Several criteria for designing unitary
constellations are derived from the analysis of the limiting
behavior. Two numerical algorithms,  SA and GA produce excellent codes
with large diversity product and diversity sum. The numerical
methods can equally be employed to optimize the diversity
function at a certain SNR.

Some algebraic structured constellations are presented. Is there
any other or better constellation structure? In~\cite{li02},
parametric codes are designed using a $A^kB^kC^k$ structure. Can
we also design other codes algebraically using other structure? It
should be an interesting question worth exploring further.

One can see that the numerical algorithms produce excellent codes
frequently although only some simple random mechanism is applied.
To explain this, we need explore further the deeper geometric
structure of $U(M)$ and the Stiefel manifold $\S_{T,M}$. Also our
implementation of the algorithm is very simple and naive. We are
still wondering if more subtle and delicate implementation will
generate more remarkable results.

We already have some examples showing that the diversity product
(or sum) doesn't necessarily mean good performance.  Basically
they are just the maximized minimum distance. We know that some of
the turbo or LDPC codes don't have large minimum distance either,
but they are capacity achieving codes which perform very well. The
reason for this is because these codes have remarkable distance
distribution. Is there any possibility to look further into the
transformation of the criterion? For instance, consider the sum of
all the distance. In this case, better distance distribution seems
to be guaranteed. In this case the problem can be converted to a
polynomial minimization problem which is widely explored in
optimization theory.

We did some elementary explorations for general type
constellations. Compared to the special form of constellations,
general form constellations don't have the unitary differential
modulation which will speed up the transmission rate. Maybe it is
possible to develop a differential modulation scheme for general
constellations? What performance could a general form
constellation have compared to a special form one? It seems like
that a lot of further research will be needed to address these
questions.

\vspace{0.3cm}
\section*{Appendix}\appendix

\Section{An Upper Bound for the  Optimal Diversity Product and
  Sum for $3$ Element Constellations} \label{three}

In this appendix, we are going to give an upper bound for the
diversity product of an arbitrary dimensional constellation with
exactly three elements. As a consequence, one can have the
optimal value and optimal conditions for a $2$ dimensional
constellation with $3$ elements. The same technique can be used
to give the optimal bound for diversity sum. The main idea
follows from~\cite{ha01p}.

Let $A=(a_{ij})_{n \times n}$ be a $n\times n$ matrix and denote
with $S_n$ the symmetric group in $n$ elements. Let
$$
F(n):=\displaystyle\sup_{A \in U(n)} \sum_{\sigma \in S_n}
|\prod_{i=1}^n a_{i\sigma(i)}|.
$$
Then one has:
\begin{thm}                                    \label{thm:1}
  Let $\V$ be a unitary constellation as a subset of $U(n)$ with
  $3$ elements, then
  \begin{equation}                     \label{Eq:bound}
     \sup_{\V \subset U(n)} \prod \V \leq \sqrt[n]{F(n)} \frac{\sqrt{3}}{2}.
\end{equation}
\end{thm}

For the proof of this theorem we will need the following
technical lemma: Let $S(m,n)$ be the set of $n \times m$ matrices
of the form $\Phi=(\Phi_{ij})$ where $0 \leq \Phi_{ij} \leq \pi$
and $\sum\limits_{i=1}^{n} \Phi_{ij} =\pi$.  Let
$d_i(\Phi):=\prod\limits_{j=1}^{m}\sin \Phi_{ij}$, then we have:
\begin{lem}                 \label{lem:1}
  $$
  \max\limits_{\Phi \in S(m,n)} \min\limits_{i=1,2,\cdots,n}
  d_i(\Phi)=\left(\sin \frac{\pi}{n}\right)^m
  $$
  and equality holds if and only if $\Phi_{ij}=\frac{\pi}{n}$,
  $i=1,2,\cdots,n$ and $j=1,2,\cdots,m$.
\end{lem}
\begin{proof} Let $d$ be the left hand side in above formula. Then
\begin{multline*}
  d^n \leq \prod_{i=1}^nd_i(\Phi)=\prod_{i=1}^n \prod_{j=1}^m
  \sin\Phi_{ij}
  = \prod_{j=1}^m \prod_{i=1}^n \sin \Phi_{ij} \\
  \leq \prod_{j=1}^m \left(\frac{\sum_{i=1}^n \sin \Phi_{ij}}{n}
  \right)^n \leq \prod_{j=1}^m \left( \sin
    \frac{\pi}{n}\right)^n= \left( \sin
    \frac{\pi}{n}\right)^{nm},
\end{multline*}
that means
$$
d \leq \left(\sin \frac{\pi}{n}\right)^m .
$$
On the other hand if one chooses $\Phi_{ij}=\frac{\pi}{n}$ one
sees that
$$
d \geq \left( \sin \frac{\pi}{n}\right)^m
$$
which establishes the claimed equality. We leave it to the
reader to verify that the maximal value is achieved in a unique
manner.
\end{proof}

\begin{proof} [Proof of Theorem~\ref{thm:1}]
  Consider a constellation $\V \subset U(n)$ with cardinality
  $|\V|=3$. Without loss of generality we can assume that $\V$
  has the form:

  $$\V = \{ I_n, D, A\},$$

  where $I_n$ is the $n\times n$ identity matrix and $D$ is a
  diagonal matrix of the form
  $$
  D=\diag(e^{i\theta_1},e^{i\theta_2},\cdots,e^{i\theta_n})
  $$
  and $A$ is an arbitrary unitary matrix.  Assume $A$ has
  eigenvalues $e^{i\varphi_1},\ldots,e^{i\varphi_n}$, i.e.  there
  is a unitary matrix $U=(u_{ij})_{n \times n}$ such that
  $$
  U^{-1}AU=\diag(e^{i\varphi_1},e^{i\varphi_2},\cdots,e^{i\varphi_n}).
  $$

  If either
  $$
  |\det(I-D)| \leq 3^\frac{n}{2} \mbox{ or } |\det(I-A)| \leq
  3^\frac{n}{2}
  $$
  then automatically we have
  $$
  \prod \V \leq \frac{\sqrt[n]{F(n)}\sqrt{3}}{2}.
  $$

  Assume therefore that $|\det(I-D)| > 3^\frac{n}{2}$ and
  $|\det(I-A)| > 3^\frac{n}{2}$, that is
  $$
  |\det(I-\diag(e^{i\theta_1},e^{i\theta_2},\cdots,e^{i\theta_n}))|
  > 3^\frac{n}{2},
  $$
  $$
  |\det(I-\diag(e^{i\varphi_1},e^{i\varphi_2},\cdots,e^{i\varphi_n})|
  > 3^\frac{n}{2},
  $$
  then according to Lemma~\ref{lem:1} we have the following
  inequality:
\begin{multline*}
  |\det(D-A)|=|\det(D-U\diag(e^{i\varphi_1},e^{i\varphi_2},
  \cdots,e^{i\varphi_n})U^{-1})|\\
  =|\det(DU-U\diag(e^{i\varphi_1},e^{i\varphi_2},
  \cdots,e^{i\varphi_n}))|
  =|\det((u_{ij}(e^{i\theta_i}-e^{i\varphi_j})))|\\
  \leq \sum_{\sigma \in S_n} \prod_{i=1}^{n}
  |u_{i\sigma(i)}(e^{i\theta_i}-e^{i\varphi_{\sigma(i)}})| \leq
  F(n) 3^\frac{n}{2}.
\end{multline*}
Taking the $n$th root and dividing the result by 2 shows also in
this case that the diversity product is at most the value on the
right hand side in\eqr{Eq:bound}.
\end{proof}

In general we do not know how sharp the estimate is.  However
when $n=2$ we have:

\begin{thm}Let $\V$ be a unitary constellation as a subset of $U(2)$ with
  $3$ elements, then
  $$
  \sup_{\V \subset U(n)} \prod \V = \sqrt{3}/2.
  $$
\end{thm}

\begin{proof}
  Since $U(2)$ has the explicit form\eqr{u2} one can verify that
  $F(2)=1.$ From Theorem~\ref{Eq:bound}, we have

  $$\displaystyle\sup_{\V \in C^2, |\V|=3} \prod \V \leq
  \sqrt{3}/2.
  $$

  On the other hand Remark~\ref{poss} shows the existence
  constellations of three elements whose diversity product is
  $\sqrt{3}/2.$
\end{proof}

\begin{rem}                                  \label{poss}
  If $\V\subset U(2)$ is a constellation of cardinality 3 having
  $\prod \V = \sqrt{3}/2$, then $\V$ must have the following
  form:
  $$\{C,CADA^{-1},CBEB^{-1}\} \mbox{ or }
  \{C,CAFA^{-1},CBGB^{-1}\} \mbox{ or }$$
  $$\{C,ADA^{-1}C,BEB^{-1}C\} \mbox{ or }
  \{C,AFA^{-1}C,BGB^{-1}C\}$$
  where
  $$
  D=\left(\begin{array}{cc}
      e^{\frac{2\pi i}{3}} & 0 \\
      0 & e^{\frac{-2\pi i}{3}} \
  \end{array}\right)   \mbox{     } E=\left(\begin{array}{cc}
    e^{\frac{4\pi i}{3}} & 0 \\
    0 & e^{\frac{-4\pi i}{3}} \
  \end{array}\right)$$
$$
F=\left(\begin{array}{cc}
    e^{\frac{2\pi i}{3}} & 0 \\
    0 & e^{\frac{4\pi i}{3}} \
  \end{array}\right) \mbox{    }
G=\left(\begin{array}{cc}
    e^{\frac{4\pi i}{3}} & 0 \\
    0 & e^{\frac{2\pi i}{3}} \
  \end{array}\right)
$$
and $A,B,C$ are arbitrary $2\times 2$ unitary matrices.
\end{rem}

\begin{rem}
  If $n=3$ then we computed with the help of Matlab $F(3) \cong
  1.299$.  This results in the upper bound:
  $$
  \sup_{\V \subset U(n)} \prod \V \leq
  \sqrt[3]{F(3)}\sqrt{3}/2\cong 0.95
  $$

  When $n \geq 4$, we believe $F(n) \geq (\frac{2}{\sqrt{3}})^n$
  and the inequality should be trivial.
\end{rem}

For the diversity sum of a $2 \times 2$ constellation with $3$
elements, we have the same result as for the diversity product.
The following Lemma deals with diagonal unitary constellation
first.

\begin{lem} \label{lem:2}
  Let $\V$ be a diagonal unitary constellation as a subset of
  $U(2)$ with $3$ elements, then we have
  $$
  \sup_{\V \subset D(n)} \sum \V = \sqrt{3}/2,
  $$

  where $D(n)$ denote the set of all the diagonal unitary
  matrices.
\end{lem}

\begin{proof}
  Without loss of generality, we assume that $\V=\{I_2,A,B\},$
  where $I_2$ is $2 \times 2$ matrix and

  $$
  A=\left(\begin{array}{cc}
      e^{i\theta_1}&0\\
      0&e^{i\theta_2}\\
          \end{array} \right), \qquad
        B=\left(\begin{array}{cc}
            e^{i\varphi_1}&0\\
            0&e^{i\varphi_2}\\
  \end{array} \right).
$$

Suppose by contradiction

$$
\sup_{\V \subset D(n)} \sum \V > \sqrt{3}/2,
$$

that is there are matrices $A, B$ such that,
\begin{eqnarray*}
  {|1-e^{i\theta_1}|}^2+{|1-e^{i\theta_2}|}^2 &> &6,\\
{|1-e^{i\varphi_1}|}^2+{|1-e^{i\varphi_2}|}^2 &> & 6,\\
{|e^{i\varphi_1}-e^{i\theta_1}|}^2+{|e^{i\varphi_2}-e^{i\theta_2}|}^2
&> & 6.
\end{eqnarray*}

One can verify that after changing of variables, the three above
equation can be rewritten in the following way:

$$
\sin^2 x + \sin^2 x_1 > \frac{3}{2},\ \sin^2 y + \sin^2 y_1 >
\frac{3}{2}\ \mbox{ and } \ \sin^2 z + \sin^2 z_1 > \frac{3}{2},
$$

where $ \frac{\pi}{4} < x,y,z,x_1,y_1,z_1 < \frac{\pi}{2}$,
$x+y+z=\pi$ and $x_1+y_1+z_1=\pi.$

By adding up three inequalities one can derive:

\begin{equation} \label{Eq:sin}
\sin^2 x + \sin^2 y + \sin^2 z + \sin^2 x_1 + \sin^2 y_1 + \sin^2
z_1 > \frac{9}{2}
\end{equation}

Since the function $\sin^2(x)$ is concave on the interval
$[\frac{\pi}{4}, \frac{\pi}{2}]$ one  shows that

\begin{multline*}
\sin^2 x + \sin^2 y + \sin^2 z + \sin^2 x_1 + \sin^2 y_1 +
\sin^2 z_1\\
\leq 3 \sin^2\frac{x+y+z}{3}+3
\sin^2\frac{x_1+y_1+z_1}{3}=\frac{9}{2},
\end{multline*}
which contradicts\eqr{Eq:sin}.

On the other hand, let

$$
u=\left\{I, A, B|A=\left(\begin{array}{cc}
      e^{\frac{2\pi i}{3}}&0\\
      0& e^{\frac{-2\pi i}{3}}\\
              \end{array}\right),
            B=\left(\begin{array}{cc}
                e^{\frac{4\pi i}{3}}&0\\
                0& e^{\frac{-4\pi i}{3}}\\
\end{array}\right)\right\},
$$

one verifies that

$$\sum u = \sqrt{3}/2.$$

That means the upper bound can be reached.

\end{proof}

With the lemma above, we can prove the following theorem.

\begin{thm} \label{thm:2} Let $\V$ be a unitary constellation as a subset of $U(2)$ with
  $3$ elements, then we have
\begin{equation} \label{newproof}
\sup_{\V \subset U(n)} \sum \V = \sqrt{3}/2.
\end{equation}
\end{thm}

\begin{proof}

  Consider a constellation $\V \subset U(2)$ with cardinality
  $|\V|=3$. Without loss of generality we can assume that $\V$
  has the form:

  $$\V = \{ I_2, D, A\},$$

  where $I_n$ is the $n\times n$ identity matrix and $D$ is a
  diagonal matrix of the form
  $$
  D=\diag(e^{i\theta_1},e^{i\theta_2})
  $$
  and $A$ is an arbitrary unitary matrix.  Assume $A$ has
  eigenvalues $e^{i\varphi_1},e^{i\varphi_2}$, i.e.  there is a
  unitary matrix $U=(u_{ij})_{2 \times 2}$ such that
  $$
  U^{-1}AU=\diag(e^{i\varphi_1},e^{i\varphi_2}).
  $$

  If either
  $$
  {\|I-D\|}_F^2 \leq 6 \mbox{ or } {\|I-A)\|}_F^2 \leq 6
  $$
  then automatically we have
  $$
  \sum \V \leq \frac{\sqrt{3}}{2}.
  $$

  Assume therefore that ${\|I-D\|}_F^2 > 6$ and ${\|I-A)\|}_F^2 >
  6$, that is
  $$
  {\|I-\diag(e^{i\theta_1},e^{i\theta_2})\|}_F^2 > 6,
  $$
  $$
  {\|I-\diag(e^{i\varphi_1},e^{i\varphi_2})\|}_F^2 > 6,
  $$
  then according to Lemma~\ref{lem:2} and the explicit
  form~\ref{u2} of $U(2)$, we have the following inequality:
\begin{multline*}%
  \|D-A\|_F^2={\|D-U\diag(e^{i\varphi_1},e^{i\varphi_2})U^{-1}\|}_F^2
  ={\|DU-U\diag(e^{i\varphi_1},e^{i\varphi_2})\|}_F^2\\
  ={\|(u_{ij}(e^{i\theta_i}-e^{i\varphi_j}))\|}_F^2
  ={|u_{11}|}^2({|e^{i\theta_1}-e^{i\varphi_1}|}^2
  +{|e^{i\theta_2}-e^{i\varphi_2}|}^2)\\
  + {|u_{22}|}^2({|e^{i\theta_1}-e^{i\varphi_2}|}^2
  +{|e^{i\theta_2}-e^{i\varphi_1}|}^2) \leq
  6({|u_{11}|}^2+{|u_{22}|}^2)=6.
\end{multline*}
Taking the $2$th root and dividing the result by $2\sqrt{2}$
shows also in this case that the diversity sum is at most the
value on the right hand side in\eqr{newproof}.

On the other hand one can take the same constellation as  in
Remark~\ref{poss} and show that the upper bound for the diversity
sum is reached.
\end{proof}


\bibliography{huge} \bibliographystyle{plain}

\end{document}